\documentclass[12pt]{amsart}
\usepackage{epsfig}
\usepackage{amssymb}

\newtheorem{theorem}{Theorem}[]

\newtheorem{corollary}[theorem]{Corollary}

\theoremstyle{definition}

\theoremstyle{remark}
\newtheorem{remark}[theorem]{Remark}

\numberwithin{equation}{section}

 \def\R{{\mathbb{R}}}
 \def\Z{{\mathbb{Z}}}
 \def\mcg{{\mathcal M}}

\begin{document}

 \title{On the linearity of certain mapping class groups}

 \author{Mustafa Korkmaz}
 \address{Department of Mathematics, Middle East Technical University,
 06531 Ankara, Turkey} \email{korkmaz@arf.math.metu.edu.tr}
 \date{\today}
 \begin{abstract}
 S. Bigelow proved that the braid groups are linear. That is, there is a
 faithful representation of the braid group into the general linear group of
 some field. Using this, we deduce from previously known results that the
 mapping class group of a sphere with punctures and hyperelliptic mapping class
 groups are linear. In particular, the mapping class group of a closed
 orientable surface of genus $2$ is linear.
 \end{abstract}

 \maketitle

 \setcounter{secnumdepth}{1}
 \setcounter{section}{0}

 \section{Introduction}

 One of the well-known open problem in the theory of mapping class
 groups is that whether these groups are linear or not (cf. \cite{bi}, Problem
 $30$, p. $220$).  A group
 is called {\em linear} if it has a faithful representation into
 $GL(n,F)$ for some field $F$ and for some integer $n$.

 Recently, S. Bigelow~\cite{b} proved that the braid groups are
 linear. The braid group $B_n$ on $n$ strings divided out by its center
 is isomorphic to a finite index subgroup of the mapping class group
 of a sphere with $n+1$ marked points. Using this, we observe that
 the mapping class group of a sphere with marked points
 and that the hyperelliptic mapping class groups, which are defined
 below, are linear. In particular, the mapping class group of a closed
 orientable surface of genus $2$ is linear. The linearity of the
 mapping class group of a surface of genus $\geq 3$ still remains open.

 \section{Preliminaries}

 We first set up the notations and state the theorems used
 in the proof of the results of
 this paper. Then we prove our results.

 Let $S$ be a compact connected orientable surface of genus $g$
 with $r$ marked points (also called {\em punctures})
 contained in the interior of $S$
 and with $s$ boundary components.
 The {\em mapping class group} $\mcg_{g,r}^s$ of $S$ is
 defined to be the group of isotopy classes of orientation
 preserving diffeomorphisms of $S$ which preserve the set of marked
 points and are the identity on the boundary. The isotopies are assumed
 to fix each marked point and each boundary point.
 We denote the group $\mcg_{g,0}$ simply by $\mcg_{g}$.

 The braid group $B_n$ on $n$ strings is the group which admits a presentation
 with generators $\sigma_1,\sigma_2,\ldots,\sigma_{n-1}$, and with the relations
 $$
 \sigma_i\sigma_j=\sigma_j\sigma_i, \mbox{ if $|i-j|\geq 2$}
 $$ and $$
 \sigma_i\sigma_{i+1}\sigma_i=\sigma_{i+1}\sigma_i\sigma_{i+1} \ .
 $$
 In fact, the group $B_n$ is isomorphic to $\mcg_{0,n}^1$,
 the mapping class group of a disc with $n$ marked points.
 The generator $\sigma_i$ is the isotopy class of a certain diffeomorphism
 of $D_n$ which interchanges $i$th and $i+1$st marked points so that its square
 is a Dehn twist.

 S. Bigelow proved the following remarkable theorem in \cite{b}.

 \begin{theorem}  \label{th:bigelow}
 The braid groups are linear.
 \end{theorem}

 For a group $G$ and for a subset $X\subseteq G$, the centralizer of
 $X$ in $G$ is defined to be
 $$
 C_G(X)=\{y\in G: xy=yx \ \mbox{ for every $x\in X$} \}.
 $$
 The center of $G$ is $C_G(G)$ and it is denoted by $C(G)$;
 $$
 C(G)=\{x\in G: xy=yx \mbox{ for every }y\in G \}.
 $$

 For a field $F$, let $F_n$ denote the space of $n\times n$
 matrices with entries in $F$. As usual, $GL(n,F)$ denotes the
 group of invertible matrices.

 \begin{theorem} [\cite{w}, Theorem 6.2.] \label{th:winf}
 Let $G$ be a subgroup of $GL(n,F)\subseteq F_n$ and $H$ a normal
 subgroup of $G$ such that $H=C_G(X)$ for some subset $X$ of
 $F_n$. Then there exists a homomorphism of $G$ into $GL(n^2,F)$
 with kernel $H$.
 \end{theorem}

 \begin{corollary}
 If $G$ is linear, then so is $G/C(G)$.
 \end{corollary}
 \begin{proof}
 Take $X=G$ in Theorem~\ref{th:winf}.
 \end{proof}

 The following theorem is probably well known to algebraists and
 can easily be proved by using the induced representation (cf.
 \cite{l}).

 \begin{theorem} \label{th:subgplintry}
 Let $G$ be a group and $H$ be a subgroup of $G$ of finite index $n$.
 Then any injective homomorphism $H\to GL(k,F)$ gives rise to
 an injective homomorphism $G\to GL(kn,F)$. In
 particular, $G$ is linear if and only if $H$ is linear.
 \end{theorem}

 \section{The results}

 We are now ready to state and prove our results of this note.

 \begin{theorem} \label{th:sphere}
 The mapping class group $\mcg_{0,n}$ of a sphere with $n$ marked points is linear
 for every $n$.
 \end{theorem}
 \begin{proof}
 If $n\leq 3$, then $\mcg_{0,n}$ is a finite group and hence it is
 linear. Hence, we assume that $n\geq 4$.

 Recall that the braid group $B_{n-1}$ is isomorphic to the
 mapping class group of a disc $D_{n-1}$ with $n-1$ marked points.
 The center of the braid group $B_{n-1}$ is the infinite
 cyclic group generated by a Dehn twist about a simple closed
 curve isotopic to the boundary component of the disc $D_{n-1}$ (cf.~\cite{bi}).
 Let us glue a disc with one marked point $x$ to the boundary
 of $D_{n-1}$ to get a sphere $S$ with $n$ marked points.
 Extending the diffeomorphisms of $D_{n-1}$ to $S$ by
 the identity gives a homomorphism $\varphi$ from  $B_{n-1}$ to
 $\mcg_{0,n}$. The image $\varphi (B_{n-1})$ of $\varphi$
 is precisely the stabilizer of $x$ under the action of
 $\mcg_{0,n}$ on the set of marked points, which is of index $n$,
 and the kernel of $\varphi$ is the center $C(B_{n-1})$ of $B_{n-1}$.

 The group $\varphi (B_{n-1})$ is isomorphic to the quotient group
 $B_{n-1}/C(B_{n-1})$. Since the group $B_{n-1}$ is linear,
 so is $\varphi (B_{n-1})$ by Theorem~\ref{th:winf}. By Theorem \ref{th:subgplintry}
 the group $\mcg_{0.n}$ is linear.
 \end{proof}

\begin{figure}[hbt]
    \begin{center}
          \epsfig{file=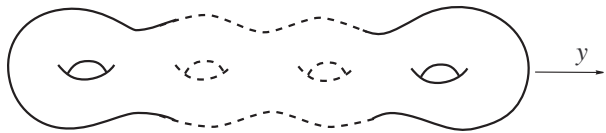,width=8.5cm}
          \caption{A surface embedded in $\R^3$ which is invariant under $J$.}
          \label{yuzey}
    \end{center}
 \end{figure}

 Suppose that a closed connected orientable surface of genus $g$ is embedded
 in the $xyz$-space as in Figure~\ref{yuzey}
 in such a way that it is invariant under the
 rotation $J(x,y,z)=(-x,y,-z)$ about the $y$-axis. Let us denote the isotopy class of
 $J$ by $\jmath$. The hyperelliptic mapping class group of genus $g$
 is defined to be the centralizer $C_{\mcg_g}(\jmath)$ of $\jmath$
 in $\mcg_g$. If $g=1$ or $2$, then the hyperelliptic mapping class
 group is equal to the mapping class group.

 \begin{theorem} \label{th:hypellp}
 Let $S$ be a closed connected orientable surface of genus $g$.
 Then the hyperelliptic mapping class group of $S$ is linear. In particular,
 the mapping class group of a closed connected orientable surface
 of genus $2$ is linear.
 \end{theorem}

 \begin{proof}
 Since the mapping class group of a torus is isomorphic to $SL(2,\Z)$,
 which is linear, we can assume that $g\geq 2$.

 There is a well known short exact sequence
 $$
 1\longrightarrow \Z_2 \longrightarrow
 C_{\mcg_g}(\jmath)\stackrel{p}\longrightarrow \mcg_{0,2g+2} \longrightarrow 1,
 $$
 where $\Z_2$ is the subgroup generated by the involution
 $\jmath$ (cf. \cite{bh}). Let $\rho:\mcg_g\to Sp(2g,3)$ the
 natural homomorphism from the mapping class group to the
 symplectic group over the finite field with three elements
 given by the action of $\mcg_g$ on the first homology group of
 $S$. Let us denote by $H$ the intersection of $C_{\mcg_g}(\jmath)$ with
 the kernel of $\rho$. Hence, $H$ is a finite index subgroup of
 $C_{\mcg_g}(\jmath)$.  As $\jmath$ acts as the minus identity on
 the first homology, it is not contained in $H$. Hence,
 the restriction of $p$ to $H$ is injective. Since $p(H)$ is
 of finite index in $\mcg_{0,2g+2}$, $H$ is linear.
 Therefore, the group $C_{\mcg_g}(\jmath)$ is linear by
 Theorem~\ref{th:subgplintry}.

 The second statement follows from the fact that
 $\mcg_2=C_{\mcg_2}(\jmath)$.
 \end{proof}

 \begin{remark}
 Bigelow proves that the braid group $B_n$ can be embedded in
 $GL(\frac{n(n-1)}{2},\Z[q^{\pm 1},t^{\pm 1}])$. Since the ring
 $\Z[q^{\pm 1},t^{\pm 1}]$ can be embedded in the field $\R$ of
 real numbers by assigning to $q,t$ two algebraically independent
 nonzero real numbers, the group $B_n$ embeds into
 $GL(\frac{n(n-1)}{2},\R)$. Using Theorems~\ref{th:winf}
 and~\ref{th:subgplintry} and the fact that the order of the group $Sp(2g,3)$
 is $3^{g^2}\prod_{i=1}^{g}(3^{2i}-1)$, it can be deduced from the proofs of
 Theorems~\ref{th:sphere} and~\ref{th:hypellp} that
 \begin{itemize}
  \item[1)] $\mcg_{0,n}$ embeds into $GL(\frac{n(n-1)^2(n-2)^2}{4},\R)$,
  \item[2)] the hyperelliptic mapping class group of genus $g$ embeds into
 $$GL(2(g+1)g^2(2g+1)^2 3^{g^2}\prod_{i=1}^{g}(3^{2i}-1),\R).$$
 \end{itemize}
 In particular, $\mcg_2$ embeds
 into $GL(2^{10}\,3^5\,5^3,\R)$.
 \end{remark}

 \noindent
 {\bf Acknowledgement.}
 The author would like to thank Mahmut Kuzucuo\u{g}lu for fruitful
 discussions.

 \bibliographystyle{amsplain}

 \end{document}